\pdfoutput=1
\documentclass[11pt]{article}
\usepackage[margin=1in]{geometry}
\usepackage{amsmath,amssymb,amsthm}
\usepackage{booktabs}
\usepackage[utf8]{inputenc}
\DeclareUnicodeCharacter{211A}{\ensuremath{\mathbb Q}}
\DeclareUnicodeCharacter{2102}{\ensuremath{\mathbb C}}
\DeclareUnicodeCharacter{2115}{\ensuremath{\mathbb N}}
\DeclareUnicodeCharacter{2124}{\ensuremath{\mathbb Z}}
\DeclareUnicodeCharacter{211D}{\ensuremath{\mathbb R}}
\DeclareUnicodeCharacter{2192}{\ensuremath{\to}}
\DeclareUnicodeCharacter{2260}{\ensuremath{\neq}}
\DeclareUnicodeCharacter{2264}{\ensuremath{\le}}
\DeclareUnicodeCharacter{2227}{\ensuremath{\wedge}}
\DeclareUnicodeCharacter{03B1}{\ensuremath{\alpha}}
\DeclareUnicodeCharacter{03C3}{\ensuremath{\sigma}}
\usepackage[hidelinks]{hyperref}

\numberwithin{equation}{section}
\newtheorem{theorem}{Theorem}[section]
\newtheorem{proposition}[theorem]{Proposition}
\newtheorem{lemma}[theorem]{Lemma}
\newtheorem{corollary}[theorem]{Corollary}
\theoremstyle{definition}
\newtheorem{definition}[theorem]{Definition}
\theoremstyle{remark}
\newtheorem{remark}[theorem]{Remark}

\newcommand{\Q}{\mathbb Q}
\newcommand{\Z}{\mathbb Z}
\newcommand{\N}{\mathbb N}

\newcommand{\C}{\mathbb C}

\newcommand{\artanh}{\operatorname{artanh}}
\newcommand{\lcm}{\operatorname{lcm}}
\newcommand{\rk}{\operatorname{rank}}
\newcommand{\red}{\operatorname{red}}
\newcommand{\lean}[1]{\texttt{\detokenize{#1}}}

\title{Logarithms of algebraic numbers:\\
the irrationality exponent and the rank of conjugate ratios}
\author{Gregory V. Chudnovsky\thanks{Developed with Claude (Anthropic). The family of
interpolation centers of Section~\ref{sec:centers} was found with Claude Fable~5.2; the
general argument, its Lean~4 formalization and this text were written with Claude
Opus~5.5. All results are formally verified in Lean~4; see Section~\ref{sec:lean}.}}
\date{November 25, 2026}

\begin{document}
\maketitle

\begin{abstract}
Let $a>0$ be a real algebraic number, $a\neq1$, with conjugates $a_1,\dots,a_d$, and let
$r(a)$ be the rank of the multiplicative group generated by the ratios $a_i/a_j$. We prove
that the irrationality exponent of $\log a$ satisfies $\mu(\log a)\le r(a)+2$. In particular
$\mu(\log a)\le d+1$, $\mu(\log a)\le3$ for every positive quadratic $a\neq1$ (for instance
$\mu(\log\frac{1+\sqrt5}2)\le3$), and $\mu(\log a)=2$ when a power of $a$ is rational. The
interpolation method of \cite{L} gives $2d$; its obstruction is that only one embedding of
the interpolation determinant is analytically small. We use centers
$\bigl(\prod_la_{i_l}^{c_l},(\sum_lc_l)\,p/q\bigr)$ indexed by a box of exponent vectors:
at \emph{every} embedding each center lies near one of a controlled number of translated
logarithmic curves, so the field norm costs no degree factor, and the number of curves is
governed by $r(a)$. All results are formally verified in Lean~4, on top of the Lean
library of \cite{L}; the final theorems depend only on the standard axioms.
\end{abstract}

\section{Introduction}

For real irrational $x$ let
\[
 \mu(x)=\sup\{\nu>0:\ |x-p/q|<q^{-\nu}\ \text{for infinitely many } p/q\in\Q\}.
\]
Multiplication by a non-zero rational number does not change $\mu$ \cite[Section~3]{L}.

Let $a>0$ be a real algebraic number, $a\ne1$, of degree $d$, with conjugates
$a=a_1,a_2,\dots,a_d$ in $\C$. Put
\[
 \Gamma'(a)=\langle a_i/a_j:\ 1\le i,j\le d\rangle\subset\C^\times,\qquad
 r(a)=\rk\Gamma'(a)=\dim_\Q\bigl(\Gamma'(a)\otimes\Q\bigr).
\]
Since $\Gamma'(a)$ is generated by $a_i/a_1$, $2\le i\le d$, we have $0\le r(a)\le d-1$.
Moreover $r(a)=0$ if and only if $a^N\in\Q$ for some $N\ge1$: if all $a_i/a_1$ are roots
of unity, then $a^N$ equals all its conjugates for a suitable $N$; the converse is clear.

\begin{theorem}\label{thm:main}
For every real algebraic $a>0$, $a\ne1$,
\[
 2\le\mu(\log a)\le r(a)+2 .
\]
\end{theorem}

\begin{corollary}\label{cor:main}
Let $a>0$, $a\ne1$, be algebraic of degree $d$.
\begin{enumerate}
\item $\mu(\log a)\le d+1$.
\item If $d=2$, then $\mu(\log a)\le3$. Equivalently, by Hilbert's Theorem~90,
$\mu(\artanh(r\sqrt D))\le3$ for rational $r\ne0$ and non-square rational $D>0$ with
$|r|\sqrt D<1$. For instance, with $\varphi=(1+\sqrt5)/2$,
\[
 \mu(\log\varphi)\le3,\qquad\mu(\log(1+\sqrt2))\le3,\qquad\mu(\log(2+\sqrt3))\le3 .
\]
\item If $a^N\in\Q$ for some $N\ge1$, then $\mu(\log a)=2$.
\item For the cubic unit $a=2\cos(2\pi/7)$, $\mu(\log a)\le4$.
\end{enumerate}
\end{corollary}

\begin{proof}
(1) $r(a)\le d-1$. (2) For $d=2$ the bound is (1). An element of norm one of
$\Q(\sqrt D)$ other than $-1$ is $(1+r\sqrt D)/(1-r\sqrt D)$ with $r\in\Q$, positive
exactly when $|r|\sqrt D<1$, and $\log\frac{1+y}{1-y}=2\artanh y$; conversely
$e^{2\artanh(r\sqrt D)}$ is such an element. Further
$\log\varphi=\artanh(1/\sqrt5)$, $\log(1+\sqrt2)=\artanh(1/\sqrt2)$,
$\log(2+\sqrt3)=2\artanh(1/\sqrt3)$. (3) $r(a)=0$. (4) $d=3$.
\end{proof}

Part~(3) is \cite[Theorem~4.1]{L} and \cite{LJZ}; here it is the case $r=0$ of the same
argument. The lower bound $\mu\ge2$ holds because $\log a$ is irrational
(Hermite--Lindemann) and Dirichlet's theorem applies; only the upper bound needs proof.

\paragraph{Relation to earlier work.}
Baker's theory of linear forms in logarithms gives effective irrationality measures of
$\log a$ for every algebraic $a$, with exponents depending on $a$ and of large size; see
Waldschmidt \cite{W78} for transcendence and approximation measures of logarithms.
Amoroso and Viola \cite{AV} obtained effective $\Q(\alpha)$-irrationality measures of
$\log\alpha$, i.e.\ approximation by elements of $\Q(\alpha)$ measured by the Weil height,
for algebraic $\alpha$ satisfying an explicit positivity condition; their examples include
real quadratic irrationals and square roots of rational numbers close to $1$. For rational
$a$ the exponent is $2$ \cite{L,LJZ}. For algebraic $a$, \cite[Theorem~4.1]{L} gives
$\mu(\log a)\le2\min_N[\Q(a^N):\Q]$, which is $2d$ when no power of $a$ has smaller degree.
The loss comes from the field norm of the interpolation determinant: only the
distinguished embedding is analytically small, while the norm counts all copies of the
arithmetic cost \cite[(6.5) and Section~9]{L}. Theorem~10.1 of \cite{L} recovers exponent
$2$ for normalized values such as $\log\alpha/\sqrt D$, where every embedding is good.
Here we make every embedding good for $\log a$ itself, by enlarging the family of
interpolation centers. Theorem~\ref{thm:main} improves $2d$ to $r(a)+2\le d+1$ for every
algebraic $a$ of degree $d\ge2$. Our bounds are uniform but not effective: the argument
assumes infinitely many very good approximations and derives a contradiction.

\paragraph{The idea.}
Fix rational approximations $r_i=p_i/q_i$ of $\lambda=\log a$. For an exponent vector $c$
consider the center
\[
 \Bigl(b_c,\ k_cr_1,\dots,k_cr_m\Bigr),\qquad b_c=\prod_la_{i_l}^{c_l},\quad k_c=\sum_lc_l ,
\]
where $a_{i_0},\dots,a_{i_r}$ are conjugates of $a$ forming a maximal independent family
in the sense of Section~\ref{sec:ranks}. An embedding $\tau$ of the Galois closure sends
some conjugate $a_{i_\tau}$ to $a$. Up to $2\pi i\Z$,
\[
 k_cr_i-\log\tau(b_c)=\sum_lc_l\bigl(\lambda-\log\tau(a_{i_l})\bigr)+k_c(r_i-\lambda),
\]
so at every embedding the center lies within $k_c|r_i-\lambda|$ of the translated curve
$X=\log Y+\kappa$, $\kappa=\sum_lc_l(\lambda-\log\tau(a_{i_l}))$. The vector $(1,a_{i_\tau})$
is a rational combination of the vectors $(1,a_{i_l})$, and the corresponding combination
of the numbers $\lambda-\log\tau(a_{i_l})$ vanishes modulo $2\pi i$. Hence the translation
forgets one direction of the box: there are about $Q^{r+1}$ centers but only about $Q^r$
translations at each embedding, when $c$ runs over a box of side $Q$. The collision
argument of \cite[Section~7]{L} survives with labels (translation, transverse grade), and
balancing the parameters allows every $\nu>r+2$ (Lemma~\ref{lem:params}).

\paragraph{Formal verification.}
Every result of this paper is formally verified in Lean~4 (Section~\ref{sec:lean}). The
formalization adds fifteen files, 3656 lines, to the Lean library of \cite{L}, uses its
interpolation theorem, truncation lemma, lcm bound and norm comparison as proved there,
and edits no existing file. The final theorems depend only on the axioms \lean{propext},
\lean{Classical.choice} and \lean{Quot.sound}. The proof below follows the formal proof.

\section{Ranks}\label{sec:ranks}

Let $L\subset\C$ be the Galois closure of $\Q(a)$, $d_L=[L:\Q]$, and $\lambda=\log a$.
For every embedding $\tau:L\to\C$ let $\gamma_\tau\in L$ be the unique element with
$\tau(\gamma_\tau)=a$; it is a conjugate of $a$, since $\tau(\gamma)=a$ forces $\gamma$
to be a root of the minimal polynomial of $a$, which splits in $L$
(\lean{exists_conj_preimage}). Every conjugate $a_i$ occurs as some $\gamma_\tau$, and
$a_i/a_j$ runs over the ratios $\gamma_\tau/\gamma_{\tau'}$.

\begin{definition}
A finite family of elements of an abelian group, written multiplicatively, is
\emph{independent} if no non-trivial product of integral powers of its members equals $1$.
We call a set $s$ of embeddings \emph{independent} if the vectors
$(1,\gamma_\tau)\in\Z\times L^\times$, $\tau\in s$, are independent, i.e.\
$\sum_{\tau\in s}n_\tau=0$ and $\prod_{\tau\in s}\gamma_\tau^{n_\tau}=1$ imply $n=0$.
\end{definition}

The rank of a finitely generated abelian group is the largest size of an independent
subfamily of any generating family. Hence $r(a)$ is the largest number of independent
ratios $\gamma_\tau/\gamma_{\tau'}$, and the largest size $\rk\Pi$ of an independent set
of embeddings is the rank of the group $\Pi\subset\Z\times L^\times$ generated by the
vectors $(1,a_i)$.

\begin{lemma}\label{lem:rank}
$\rk\Pi\le r(a)+1$.
\end{lemma}

\begin{proof}
Let $s$ be independent and non-empty, fix $\tau_0\in s$, and suppose
$\prod_{\tau\in s\setminus\{\tau_0\}}(\gamma_\tau/\gamma_{\tau_0})^{n_\tau}=1$. Put
$n_{\tau_0}=-\sum_{\tau\ne\tau_0}n_\tau$. Then $\sum_{\tau\in s}n_\tau=0$ and
$\prod_{\tau\in s}\gamma_\tau^{n_\tau}=1$, so $n=0$. Hence the $|s|-1$ ratios
$\gamma_\tau/\gamma_{\tau_0}$ are independent, and $|s|-1\le r(a)$.
\end{proof}

\begin{lemma}\label{lem:relation}
Let $s$ be an independent set of embeddings of the largest size. For every embedding
$\sigma$ there are $a_\sigma\in\Z\setminus\{0\}$ and integers $n_\tau$, $\tau\in s$, with
$\sum_{\tau\in s}n_\tau=a_\sigma$ and $\prod_{\tau\in s}\gamma_\tau^{n_\tau}=\gamma_\sigma^{a_\sigma}$.
Consequently, for every embedding $\sigma$,
\begin{equation}\label{eq:relation}
 \prod_{\tau\in s}\Bigl(\frac{e^\lambda}{\sigma(\gamma_\tau)}\Bigr)^{n_\tau}
 =\frac{e^{a_\sigma\lambda}}{\sigma(\gamma_\sigma)^{a_\sigma}}=1,\qquad(n_\tau)\ne0 .
\end{equation}
\end{lemma}

\begin{proof}
If $\sigma\in s$ take $a_\sigma=1$ and $n=\mathbf 1_\sigma$. Otherwise $s\cup\{\sigma\}$
is not independent: there is $n\ne0$ on $s\cup\{\sigma\}$ with $\sum n=0$ and
$\prod\gamma^n=1$; $n_\sigma\ne0$ because $s$ is independent, and $a_\sigma=-n_\sigma$
works. Finally $\sigma(\gamma_\sigma)=a=e^\lambda$.
\end{proof}

\section{Inputs from \texorpdfstring{\cite{L}}{[L]}}\label{sec:inputs}

Fix $m\ge1$, weights $W=(w_0,w_1,\dots,w_m)\in\Q_{>0}^{m+1}$, $\theta\in(0,1)$, $v_0>0$,
and put $v_i=w_i/\theta$ for $i\ge1$. Write $w\cdot\gamma=\sum_{i\ge1}w_i\gamma_i$ and
$w_*=\min_{i\ge1}w_i$. Let $\mathcal P_W(H)$ be the space of polynomials with monomial
basis $Y^hX^\gamma$, $w_0h+w\cdot\gamma\le H$, and let
$\mathcal J_V(H)=\C[[t,u_1,\dots,u_m]]/F_V^{\ge H}$ be the jets with monomials
$t^su^\beta$ of weight $v_0s+\sum v_i\beta_i<H$.

\begin{theorem}[{\cite[Theorem~5.1]{L}}]\label{thm:interp}
Fix integers $m,Z,S\ge1$ with $S\le Z$, and positive rationals $w_0,v_0,\theta$ with
\begin{equation}\label{eq:volume}
 \theta<1,\qquad S\theta^m<1,\qquad Z\frac{w_0}{v_0}\theta^m<1.
\end{equation}
There are successive lower thresholds for positive rational weights $w_1,\dots,w_m$, each
depending only on these data and the earlier weights, with the following property. For
every admissible choice of weights there is an integer $R_0\ge1$, depending only on the
weights, such that for every collection of centers
$(b_\zeta,c_{\zeta1},\dots,c_{\zeta m})\in\C^\times\times\C^m$, $\zeta=1,\dots,Z$, in which
at most $S$ of the $b_\zeta$ share any one value and, within each group with the same
$b_\zeta$, the $c_{\zeta i}$ are pairwise distinct for every $i$, there is $H_0$ such that
\[
 \mathcal P_W(H)\to\bigoplus_{\zeta}\mathcal J_V(H),\qquad
 P\mapsto\bigl[P\bigl(b_\zeta(1+t),(c_{\zeta i}+u_i+\log(1+t))_i\bigr)\bigr]_\zeta
\]
is surjective for every integer $H\ge H_0$ divisible by $R_0$.
\end{theorem}

\begin{lemma}[{Truncation, \cite[Lemma~6.1]{L}}]\label{lem:shear}
Let $T_i\ge1$ be integers with $v_0T_i\ge v_i$ and let
$G_i(t)=\sum_{1\le l<T_i}(-1)^{l+1}t^l/l$. Fix centers, take as rows all target
coefficients $t^su^\beta$ of all centers, and fix any selection of as many columns as
there are rows. Then the square coefficient matrix built with $c_{\zeta i}+G_i(t)+u_i$
and the one built with $c_{\zeta i}+\log(1+t)+u_i$ have the same determinant.
\end{lemma}

We also use $\log\lcm(1,\dots,n)\le\Lambda n$ for all $n\ge1$, with $\Lambda=\log4+4$.

\section{The centers}\label{sec:centers}

The argument needs, for every $Q\ge1$, a family of centers with the following
properties. All constants $c_K,c_E,c_N,c_R$ are positive integers independent of $Q$.

\begin{definition}\label{def:supply}
A \emph{center supply of rank $r$} for $\lambda$ in a number field $L$, with denominator
$D\ge1$, is the assignment to every $Q\ge1$ of distinct pairs $(k_z,b_z)\in\N\times
L^\times$, $z=1,\dots,Z$, an integer $e\ge0$, and for every embedding $\tau$ complex
numbers $\ell_\tau(z)$, such that
\begin{enumerate}
\item[(S1)] $Z\ge Q^{r+1}$;
\item[(S2)] $k_z\le c_KQ$, and $e\le c_EQ$ with $D^eb_z$ integral for all $z$;
\item[(S3)] $e^{\ell_\tau(z)}=\tau(b_z)$ and $|\ell_\tau(z)|\le c_RQ$;
\item[(S4)] the translations $\kappa_\tau(z)=k_z\lambda-\ell_\tau(z)$ satisfy
$|\kappa_\tau(z)|\le c_RQ$ and take at most $c_NQ^r$ values for each $\tau$.
\end{enumerate}
\end{definition}

\begin{proposition}\label{prop:supply}
Let $a$, $L$, $\lambda$ be as in Section~\ref{sec:ranks}, $D\ge1$ with $Da$ integral, and
$s=\{\tau_0,\dots,\tau_\rho\}$ an independent set of embeddings of the largest size,
$\rho+1=\rk\Pi$. Put $\beta_l=\gamma_{\tau_l}$. For $c\in\{0,\dots,Q-1\}^{\rho+1}$ let
\[
 b_c=\prod_{l=0}^\rho\beta_l^{c_l},\qquad k_c=\sum_lc_l,\qquad
 \kappa_\tau(c)=\red\Bigl(\sum_lc_lw_{\tau,l}\Bigr),\qquad \ell_\tau(c)=k_c\lambda-\kappa_\tau(c),
\]
where $w_{\tau,l}$ is a fixed logarithm of $e^\lambda/\tau(\beta_l)$ and
$\red(x)=x-2\pi i\lfloor\operatorname{Im}x/2\pi\rfloor$. This is a center supply of rank
$\rho$, with $Z=Q^{\rho+1}$, $e=(\rho+1)Q$, $c_K=c_E=\rho+1$.
\end{proposition}

\begin{proof}
\emph{Distinctness and (S1).} If $(k_c,b_c)=(k_{c'},b_{c'})$ then $n=c-c'$ satisfies
$\sum n_l=0$ and $\prod\beta_l^{n_l}=1$, so $c=c'$ by independence. Hence $Z=Q^{\rho+1}$.

\emph{(S2).} $k_c\le(\rho+1)Q$. Since $D\gamma_\tau$ is integral for every $\tau$ (it is a
conjugate of $Da$) and $\sum_lc_l\le e$,
$D^eb_c=D^{e-\sum c_l}\prod_l(D\beta_l)^{c_l}$ is integral.

\emph{(S3).} $\red(x)\equiv x\bmod 2\pi i\Z$, so
$e^{\ell_\tau(c)}=e^{k_c\lambda}\prod_l\bigl(e^\lambda/\tau(\beta_l)\bigr)^{-c_l}=\prod_l\tau(\beta_l)^{c_l}=\tau(b_c)$.
With $W=\sum_{\tau,l}|w_{\tau,l}|$ we have $|\red(x)|\le|x|+2\pi$, hence
$|\kappa_\tau(c)|\le QW+2\pi$ and $|\ell_\tau(c)|\le(\rho+1)Q|\lambda|+QW+2\pi$; take
$c_R=\lceil(\rho+1)|\lambda|+W+8\rceil$.

\emph{(S4).} Fix $\tau$. By Lemma~\ref{lem:relation} there is $n\in\Z^{\rho+1}$, $n\ne0$,
with $\prod_l(e^\lambda/\tau(\beta_l))^{n_l}=1$; replacing $n$ by $-n$ we may assume
$n_{l_0}>0$ for some $l_0$. Then $\sum_ln_lw_{\tau,l}=2\pi it$ with $t\in\Z$. Let $N_b$ be
an integer bounding all $|n_l|$, for all $\tau$. For $c$ in the box let
$s_0=\lfloor c_{l_0}/n_{l_0}\rfloor\in[0,Q)$ and $c'=c-s_0n\in\Z^{\rho+1}$. Then
$\sum c'_lw_{\tau,l}=\sum c_lw_{\tau,l}-2\pi is_0t$, so
$\red(\sum c'_lw_{\tau,l})=\kappa_\tau(c)$, while $0\le c'_{l_0}<n_{l_0}\le N_b$ and
$|c'_l|\le Q+QN_b$ for $l\ne l_0$. Hence
\[
 \#\{\kappa_\tau(c)\}\le N_b\,(Q(2N_b+1)+1)^\rho\le N_b(2N_b+2)^\rho\,Q^\rho .\qedhere
\]
\end{proof}

By Lemma~\ref{lem:rank}, $\rho\le r(a)$. Theorem~\ref{thm:main} therefore follows from
Proposition~\ref{prop:supply} and the following theorem, applied with $r=\rho$.

\begin{theorem}\label{thm:grid}
If $\lambda\ne0$ is real and has a center supply of rank $r$ in some number field, then
$\mu(\lambda)\le r+2$.
\end{theorem}

The rest of the paper proves Theorem~\ref{thm:grid}.

\begin{remark}
For a quadratic $\alpha>0$ of norm $\pm1$, $r=1$ and the centers are, up to a change of
indices, the grid $(\alpha^j,k\,p_i/q_i)$ with translations $(k-j)\lambda$ at the identity
and $(k+j)\lambda$ at the conjugate embedding (after replacing $\alpha$ by $\alpha^2$ if
the norm is $-1$). This case was formalized first, directly
(\lean{quadratic_norm_pm_one_log_irrationalityExponent_bounds}).
\end{remark}

\section{The construction}\label{sec:construction}

Let $\lambda$ have a center supply of rank $r$ in $L$ with denominator $D$. Fix $\nu>r+2$
and suppose, to reach a contradiction, that for every $Q_0$ there is $p/q$ with $q\ge Q_0$
and
\begin{equation}\label{eq:approx}
 |\lambda-p/q|\le q^{-\nu}.
\end{equation}

\paragraph{Fixed parameters.}
Choose rationals $\theta,A,B,C$ as in Lemma~\ref{lem:params}:
\begin{equation}\label{eq:P}
 0<\theta<A<1,\quad \nu(A-\theta)>1-\theta,\quad
 AC^{r}<B,\quad CB<1,\quad B<C^{r+1}\theta,\quad C\theta<1 .
\end{equation}
Then $C>A/\theta>1$. Let $g_0=\nu(A-\theta)-(1-\theta)>0$,
$\eta=\min\bigl(\tfrac12,g_0/(2\nu A)\bigr)$ and $\epsilon=\min(g_0/2,1/2)$; then
\begin{equation}\label{eq:gap0}
 g:=\nu\bigl(A(1-\eta)-\theta\bigr)-(1-\theta)\ \ge\ g_0/2\ \ge\ \epsilon .
\end{equation}

\paragraph{Dimension.}
For an integer $m\ge1$, fixed in Section~\ref{sec:end}, put $Q=\lfloor C^m\rfloor$ and
take the centers $(k_z,b_z)_{z\le Z}$ of the supply for this $Q$. Put
\begin{equation}\label{eq:m}
 w_0=B^{-m},\quad v_0=2Z\theta^mw_0,\quad S=\min(c_KQ+1,Z),\quad N=(c_N+1)Q^r,\quad
 R=2(c_RQ+1).
\end{equation}
Then $Z(w_0/v_0)\theta^m=1/2$. Centers with the same $b_z$ have distinct $k_z\in[0,c_KQ]$,
so at most $S$ of them share a value of $b_z$.

\paragraph{Approximations.}
With $m$ fixed, choose fractions $r_i=p_i/q_i$ satisfying \eqref{eq:approx}, $q_i\ge2$
and $p_i\ne0$, successively for $i=1,\dots,m$, so that $w_i:=\lceil\log q_i\rceil$
exceeds the $i$-th threshold of Theorem~\ref{thm:interp} and a lower bound $X+1$ fixed in
Section~\ref{sec:end}. Put
\[
 \varepsilon_i=r_i-\lambda,\quad |\varepsilon_i|\le2q_i^{-\nu},\quad
 T_i=\lceil w_i/(\theta v_0)\rceil,\quad
 G_i(t)=\sum_{1\le l<T_i}\frac{(-1)^{l+1}t^l}{l},\quad L_i=\lcm(1,\dots,T_i-1).
\]

\paragraph{Matrix and minor.}
For each $z$ take the center $(b_z,k_zr_1,\dots,k_zr_m)$. Within a group with the same
$b_z$ the $k_z$ are distinct, and $r_i\ne0$, so the $k_zr_i$ are distinct; by
Theorem~\ref{thm:interp} (the volume conditions hold once $S\theta^m<1$) there are $R_0$
and $H_0$. Let $H\ge H_0$ run through multiples of $R_0$. Columns: $Y^hX^\gamma$ with
$w_0h+w\cdot\gamma\le H$. Rows: $(z,s,\beta)$ with $v_0s+w\cdot\beta/\theta<H$. Entry:
\begin{equation}\label{eq:entry}
 [t^su^\beta]\ b_z^h(1+t)^h\prod_{i=1}^m\bigl(k_zr_i+G_i(t)+u_i\bigr)^{\gamma_i}\ \in L .
\end{equation}
As $v_0T_i\ge v_i$, Theorem~\ref{thm:interp} and Lemma~\ref{lem:shear} give a non-zero
maximal minor $\Delta=\Delta_H\in L$ using all $M=M_H$ rows. Put
$\bar b=\frac1{MH}\sum_{\rm rows}w\cdot\beta\in[0,\theta)$. As $H\to\infty$ with
everything else fixed,
\begin{equation}\label{eq:counts}
 M\sim\frac{Z\theta^mH^{m+1}}{(m+1)!\,v_0\prod_iw_i},\qquad
 J_H:=\#\{\delta:w\cdot\delta\le AH\}\sim\frac{A^mH^m}{m!\prod_iw_i}.
\end{equation}

\section{Arithmetic lower bound}\label{sec:arith}

\begin{lemma}\label{lem:arith}
\begin{equation}\label{eq:lower}
 \frac1{d_L}\sum_\tau\frac{\log|\tau\Delta|}{MH}\ \ge\ -(1-\bar b)-E_{\rm ar},\qquad
 E_{\rm ar}=\frac{\Lambda m}{\theta v_0}+\Lambda\sum_{i=1}^m\frac1{w_i}+\frac{\theta}{w_*}
 +\frac{e\log D}{w_0},
\end{equation}
the sum over all $d_L$ embeddings of $L$.
\end{lemma}

\begin{proof}
Scale column $(h,\gamma)$ by $D^{eh}\prod_i(q_iL_i)^{\gamma_i}$ and row $(z,s,\beta)$ by
$\prod_i(q_iL_i)^{-\beta_i}$. By (S2) $D^{eh}b_z^h$ is integral, and
$q_iL_i(k_zr_i+G_i(t))=k_zp_iL_i+q_iL_iG_i(t)\in\Z[t]$. Hence $C_H\Delta$ is a non-zero
algebraic integer, where $C_H=D^{e\sum_{\rm col}h}\prod_i(q_iL_i)^{J_i}$ and
$J_i=\sum_{\rm col}\gamma_i-\sum_{\rm rows}\beta_i$. Its norm is a non-zero integer, so
$\sum_\tau\log|\tau\Delta|\ge-d_L\log C_H$. Finally $\sum_{\rm col}h\le MH/w_0$,
$\sum_iJ_i\log L_i\le MH\Lambda(m/(\theta v_0)+\sum_i1/w_i)$ by the lcm bound, and
$\sum_iJ_i\log q_i\le MH(1-\bar b+\theta/w_*)$, as in \cite[Section~6]{L}.
\end{proof}

\section{Analytic upper bound at every embedding}\label{sec:analytic}

For a polynomial $P$, $\kappa\in\C$ and $\delta\in\Z^m_{\ge0}$ put
\[
 f^{\kappa}_{\delta,P}(z)=[u^\delta]\,P\bigl(e^z,(z+\kappa+u_i)_i\bigr),\qquad
 f^{\kappa}_{\delta,Y^hX^\gamma}(z)=\binom\gamma\delta e^{hz}\prod_i(z+\kappa)^{\gamma_i-\delta_i}
\]
($0$ unless $\delta\le\gamma$).

\begin{lemma}[Row identity]\label{lem:row}
Let $\tau$ be an embedding, $\ell=\ell_\tau(z)$ and $\kappa=\kappa_\tau(z)$. For every
polynomial $P$,
\[
 [t^su^\beta]\,P\bigl(\tau(b_z)(1+t),(k_zr_i+\log(1+t)+u_i)_i\bigr)
 =\sum_{\delta\ge\beta}\binom\delta\beta\prod_i(k_z\varepsilon_i)^{\delta_i-\beta_i}\,
 [t^s]\,f^{\kappa}_{\delta,P}\bigl(\ell+\log(1+t)\bigr).
\]
\end{lemma}

\begin{proof}
With $z=\ell+\log(1+t)$ we have $e^z=\tau(b_z)(1+t)$ by (S3), and
$k_zr_i+\log(1+t)=z+(k_z\lambda-\ell)+k_z\varepsilon_i=z+\kappa+k_z\varepsilon_i$. Expand
$P$ in $u$.
\end{proof}

Since $k_zr_i\in\Q$ is fixed by $\tau$, the matrix of $\tau\Delta$, computed with the
exact logarithm (Lemma~\ref{lem:shear}), has exactly these entries: every embedding is
good.

\begin{proposition}\label{prop:upper}
For every embedding $\tau$ and all large $H$ divisible by $R_0$,
\[
 \frac{\log|\tau\Delta|}{MH}\le E_{\rm an}+o(1)
 +\max\Bigl\{-c\,\frac{\eta^2M}{HJ_HN},\ -\nu\bigl(A(1-\eta)-\bar b\bigr)\Bigr\},
\]
where $c=\frac{\log2}{4}$, the term $o(1)$ does not depend on $\tau$, and
\[
 E_{\rm an}=\frac{\log2+\log(2(c_KQ+1))+\nu}{w_*}+\frac R{w_0}+\frac{\log(4R)}{w_*}
 +\frac{\log2}{v_0}.
\]
\end{proposition}

\begin{proof}
Expand the rows by Lemma~\ref{lem:row}; the factors $(k_z\varepsilon_i)^{\delta_i-\beta_i}$
are controlled by $k_z\le c_KQ$ and $|\varepsilon_i|\le2q_i^{-\nu}$. The determinant becomes
a sum over choices of $\delta$ for each row. In the approximation alternative,
$\sum w\cdot(\delta_r-\beta_r)\ge MH(A(1-\eta)-\bar b)$, which gives the second term. In the
collision alternative, at least $\eta M$ rows have $w\cdot\delta\le AH$; their labels
$(\kappa,\delta)$ take at most $NJ_H$ values by (S4), and rows with equal labels and equal
$\ell$ have equal coefficient vectors. The collision estimate of \cite[Lemma~7.1]{L},
applied with these labels on $|z|=R$ (where $|\kappa|,|\ell|\le R/2$ by (S3), (S4)) and
on $|t|=1/2$, together with Cauchy--Schwarz, gives the first term. The holomorphic and
counting factors give $E_{\rm an}$.
\end{proof}

By \eqref{eq:counts},
\begin{equation}\label{eq:Lprime}
 c\,\frac{\eta^2M}{HJ_HN}\ \longrightarrow\ L'=\frac{c\,\eta^2Z\theta^m}{(m+1)v_0A^mN}
 =\frac{c\,\eta^2B^m}{2(m+1)A^m(c_N+1)Q^r}\qquad(H\to\infty).
\end{equation}

\section{Parameters and the end of the proof}\label{sec:end}

\begin{lemma}\label{lem:params}
For every $r\ge0$ and $\nu>r+2$ there are rationals $\theta,A,B,C$ with \eqref{eq:P}.
\end{lemma}

\begin{proof}
Choose a rational $b$ with $\frac{r+1}{r+2}<b<1-\frac1\nu$ and a rational $\delta$ with
$0<\delta<\min\bigl(1,\frac{(r+2)b-(r+1)}{(r+1)(r+2)b}\bigr)$, and put $\theta=1-\delta$,
$A=1-b\delta$. Then $\theta<A<1$ and $\nu(A-\theta)=\nu(1-b)\delta>\delta=1-\theta$. By
Bernoulli's inequality, $\theta^{r+1}\ge1-(r+1)\delta$ and
$A^{r+2}(1+(r+2)b\delta)\le A^{r+2}(1+b\delta)^{r+2}=(1-b^2\delta^2)^{r+2}\le1$; since
$(1-(r+1)\delta)(1+(r+2)b\delta)=1+\delta\bigl((r+2)b-(r+1)-(r+1)(r+2)b\delta\bigr)>1$, we
get $A^{r+2}<\theta^{r+1}$. Thus $C_0=A/\theta$ satisfies $AC_0^{r+1}<1$ and $C_0\theta<1$,
and by continuity a rational $C$ slightly larger than $C_0$ satisfies $AC^{r+1}<1$,
$C\theta<1$ and $C>A/\theta>1$. Then $AC^r<1/C$ and $AC^r<C\theta\,C^r=C^{r+1}\theta$; choose
a rational $B$ with $AC^r<B<\min(1/C,C^{r+1}\theta)$.
\end{proof}

\begin{proof}[Proof of Theorem~\ref{thm:grid}]
Average Proposition~\ref{prop:upper} over the $d_L$ embeddings and compare with
\eqref{eq:lower}:
\[
 \min\Bigl\{c\,\frac{\eta^2M}{HJ_HN},\ \nu\bigl(A(1-\eta)-\bar b\bigr)\Bigr\}
 \le 1-\bar b+E_{\rm ar}+E_{\rm an}+o(1).
\]
As $\nu>1$ and $\bar b<\theta$, $\nu(A(1-\eta)-\bar b)-(1-\bar b)\ge g$. Letting
$H\to\infty$ and using \eqref{eq:Lprime}, it suffices to arrange
\begin{equation}\label{eq:gap}
 E_{\rm ar}+E_{\rm an}<\epsilon\le g,\qquad L'>1+\epsilon .
\end{equation}
Split $E_{\rm ar}+E_{\rm an}=E_1+E_2$ with
\[
 E_1=\frac{\Lambda m/\theta+\log2}{v_0}+\frac{e\log D+R}{w_0},\qquad
 E_2=\Lambda\sum_i\frac1{w_i}+\frac{\theta+\log2+\log(2(c_KQ+1))+\nu+\log(4R)}{w_*}.
\]
\emph{Choice of $m$.} By (S1), $v_0\ge2Q^{r+1}\theta^mB^{-m}\ge2^{-r}(C^{r+1}\theta/B)^m$
since $Q\ge C^m/2$, and $C^{r+1}\theta/B>1$; so the first term of $E_1$ tends to $0$. By
(S2), $e\log D+R\le(c_E\log D+2c_R+2)C^m$, so the second term is at most
$(c_E\log D+2c_R+2)(CB)^m\to0$. By \eqref{eq:Lprime} and $Q\le C^m$,
\[
 L'\ge\frac{c\,\eta^2}{2(m+1)(c_N+1)}\Bigl(\frac B{AC^r}\Bigr)^m\to\infty .
\]
Finally $S\theta^m\le(c_K+1)(C\theta)^m\to0$. Fix $m$ with $C^m\ge2$, $E_1<\epsilon/2$,
$L'>2$ and $S\theta^m<1$.

\emph{Choice of the approximations.} With $m$, $Q$ and the centers fixed, the numerator
of the last term of $E_2$ is a constant; choose $X$ with
$\Lambda\sum_i1/w_i+(\text{that constant})/w_*<\epsilon/2$ whenever
$w_*>X$, and choose the $r_i$ with $w_i>X+1$ and above the thresholds. Then
$E_{\rm ar}+E_{\rm an}<\epsilon$, which with $\epsilon\le\frac12$ gives \eqref{eq:gap}.
This contradicts the existence of the minors $\Delta_H$ for large $H$. Hence
$\mu(\lambda)\le r+2$. The order of choices is
\[
 \nu;\ (\theta,A,B,C,\eta,\epsilon);\ m,Q;\ \text{centers},\ w_0,v_0;\ \text{thresholds};\ X;\
 (r_i,w_i)_{i\le m};\ R_0,H_0;\ H\to\infty .\qedhere
\]
\end{proof}

\section{Formal verification}\label{sec:lean}

\paragraph{Setting.}
The formalization extends the Lean library of \cite{L} (Lean~4.34.1, Mathlib as recorded
there) by fifteen new files \lean{QuadraticPi/LogGrid*.lean} in the namespace
\lean{QuadraticPi.LogGrid}; no existing file is modified. The new files accompany this
paper. It uses from the library the bounded-fiber interpolation theorem
(Theorem~\ref{thm:interp})
\begin{center}\small\lean{MovingLogOccupancyInterpolation.bounded_fiber_interpolation},\end{center}
the truncation and jet-transfer lemmas (Lemma~\ref{lem:shear}), the lcm bound
(\lean{Arithmetic.log_lcmUpto_le}), the norm comparison over a number field,
the cofinal comparison of determinant bounds, the weight selection and the existence of
approximations above thresholds. The exponent is expressed by
\texttt{EventualLowerBoundAbove x $\kappa$}: for every $\nu>\kappa$ there is $Q\ge2$ with
$q^{-\nu}\le|x-p/q|$ for all integers $p$ and all $q\ge Q$; the library derives
$2\le\mu(x)\le\kappa$ from it.

\paragraph{Main statements.}
In Lean the Galois closure is an arbitrary Galois number field $E$, the number $a$ is the
image $\sigma_0(\alpha)$ of some $\alpha\in E$ under an embedding $\sigma_0:E\to\C$, and
$\sigma_0(\alpha)=e^\lambda$ with $\lambda$ real. The conjugate $\gamma_\tau$ is
\texttt{conjOf $\alpha$ $\sigma_0$ $\tau$}, and $r(a)$ is \texttt{ratioRank (conjOf $\alpha$ $\sigma_0$)}, the largest number of
independent ratios $\gamma_\tau/\gamma_{\tau'}$.
\begin{small}
\begin{verbatim}
theorem log_galois_ratioRank_irrationalityExponent_bounds [IsGalois ℚ E]
    (α : E) (σ0 : E →+* ℂ) (lam : ℝ) (hlam : lam ≠ 0)
    (hα : σ0 α = Complex.exp (lam : ℂ))
    (D : ℕ) (hD : 0 < D) (hint : IsIntegral ℤ ((D : E) * α)) :
    2 ≤ irrationalityExponent lam ∧
      irrationalityExponent lam ≤ (ratioRank (conjOf α σ0) : ℝ) + 2

theorem log_galois_degree_irrationalityExponent_bounds [IsGalois ℚ E] ... :
    2 ≤ irrationalityExponent lam ∧
      irrationalityExponent lam ≤ ((minpoly ℚ α).natDegree : ℝ) + 1
\end{verbatim}
\end{small}
The first is Theorem~\ref{thm:main}; the second is Corollary~\ref{cor:main}(1), and with
$\deg\alpha=2$ also Corollary~\ref{cor:main}(2).

\paragraph{Correspondence.}
\begin{center}\footnotesize
\begin{tabular}{@{}ll@{}}
\toprule
This paper & Lean declaration (file)\\
\midrule
Theorem~\ref{thm:main} & \lean{log_galois_ratioRank_irrationalityExponent_bounds} (\lean{LogGridRatioRank})\\
Corollary~\ref{cor:main}(1) & \lean{log_galois_degree_irrationalityExponent_bounds} (\lean{LogGridRatioRank})\\
Lemma~\ref{lem:rank} & \lean{conjRank_le_ratioRank} (\lean{LogGridRatioRank})\\
Lemma~\ref{lem:relation} & \lean{exists_max_conjIndep}, \lean{conj_span} (\lean{LogGridRank})\\
$\gamma_\tau$ & \lean{exists_conj_preimage}, \lean{conjOf}, \lean{isIntegral_conjOf} (\lean{LogGridGalois})\\
Proposition~\ref{prop:supply} & \lean{centerSupply_of_basis}, \lean{card_image_red_le} (\lean{LogGridSupply})\\
Definition~\ref{def:supply} & \lean{CenterSupply} (\lean{LogGridGeneral})\\
Theorem~\ref{thm:grid} & \lean{general_grid_irrationalityExponent_bounds} (\lean{LogGridGeneral})\\
Lemma~\ref{lem:params} & \lean{exists_grid_parameters_r}, \lean{bernoulli_gap} (\lean{LogGridParametersR})\\
choice of $m$ & \lean{exists_grid_dimension_r} (\lean{LogGridParametersR})\\
Lemma~\ref{lem:arith} & \lean{gridSelectedMinor_average_lower_bound} (\lean{LogGridNumberField})\\
Lemma~\ref{lem:row}, Proposition~\ref{prop:upper} & \lean{grid_minor_analytic_bound} (\lean{LogGridAnalytic})\\
non-zero minors & \lean{cofinal_grid_minors} (\lean{LogGridMinors})\\
end of the proof & \lean{grid_family_contradiction} (\lean{LogGridAssembly})\\
quadratic, norm $\pm1$ & \lean{quadratic_norm_pm_one_log_irrationalityExponent_bounds} (\lean{LogGridTheorem})\\
\bottomrule
\end{tabular}
\end{center}

\paragraph{Differences in presentation.}
The formal proof keeps the analytic estimate in the form of a remainder $o(1)$ that tends
to $0$ along the degrees $H$, and compares the two sides through the library's cofinal
comparison lemma. The bound $Q\ge C^m/2$, the reduction modulo $2\pi i$ and the counting
of translations are proved exactly as in Section~\ref{sec:centers}. The constants are
those of the text: $\Lambda=\log4+4$, $c=\frac{\log2}4$, $N=(c_N+1)Q^r$, $R=2(c_RQ+1)$.

\paragraph{Audit.}
The axioms used by the final theorems, as reported by \texttt{\#print axioms}, are
\lean{propext}, \lean{Classical.choice} and \lean{Quot.sound}, for
Theorem~\ref{thm:main}, Corollary~\ref{cor:main}(1), Theorem~\ref{thm:grid}, the rank-one
theorem \lean{logGrid_irrationalityExponent_bounds} and the quadratic theorem. The new
files contain no \lean{sorry}, \lean{admit}, \lean{native_decide} or axiom declarations.
Their SHA-256 checksums are:
\begin{center}\scriptsize
\begin{tabular}{@{}lll@{}}
\toprule
File & Lines & SHA-256\\
\midrule
\lean{LogGridMatrix} & 55 & \lean{a9e4615d4b41c5b19b8227c522fd372dee47ed88ead4d0eb05b774a3ceae08f8}\\
\lean{LogGridJetTransfer} & 148 & \lean{39823c125be494e7bc6cf0de6bd6d63328bc5ffab316c86067bea75875941d6d}\\
\lean{LogGridMinors} & 159 & \lean{e3393aa6aac55feef49032fff6c051c6939fa8eb824685a2bd812b277795c845}\\
\lean{LogGridShiftedColumns} & 131 & \lean{7798fd3373d3ffb995c864e0ea62f55bf4a2ee9f2e6588e68b8ca88182621cbf}\\
\lean{LogGridAnalytic} & 498 & \lean{0a2a331d622d24f7712ddccd480b18ad9d6ee17bf32e5e01ac5d5d9d5783ebea}\\
\lean{LogGridNumberField} & 242 & \lean{4ab0ba6310260a3cccfbd1d5c9b4ce60174b568f1adffe8b629425946b0247fc}\\
\lean{LogGridAssembly} & 253 & \lean{8b2a18c31aa28569894aed7dc78cb5bd1042c5889b4c4f5bb912f4b65a48649d}\\
\lean{LogGridParameters} & 225 & \lean{f290d636d5e3a0dbe3a67315a6d8038341e100aced962b22d1ed85f120693e96}\\
\lean{LogGridTheorem} & 504 & \lean{bd5c1072a5094b4c8e76316c699e10ae2128510c386ebbe0db508fbc4d2240cd}\\
\lean{LogGridParametersR} & 319 & \lean{45de44ca81f75a766e3ee22fb8f25ac6dcfdb72b4f4c1d1a55f6ec6ed06b5f42}\\
\lean{LogGridGeneral} & 317 & \lean{e5a1144ee496068343df9ea5720918f2e9ee225433182d7cefcd2db032d23158}\\
\lean{LogGridSupply} & 356 & \lean{ee500f12828082fff95def16f491254a1c622d13257373cb54f9b8204d2ee144}\\
\lean{LogGridRank} & 174 & \lean{53f726e704d5cd3fd54deb98ee2e81b05c3e399c1688e76c4a456877ce5f6d14}\\
\lean{LogGridGalois} & 145 & \lean{1f0d37f556da8aba3b32bba8bcc7c42dcc835b9cb38166a0725bec1449d838cd}\\
\lean{LogGridRatioRank} & 130 & \lean{c5846eba58a79372f4d2e2a53b659bb0fbbaae3e8c95cc1de1529315ec1d6be0}\\
\bottomrule
\end{tabular}
\end{center}

\section{Remarks}

\begin{remark}[Special cases]
For $r(a)=0$ the family is $\{(a^k,kr_i)\}$ up to roots of unity, the diagonal family of
\cite{L} and \cite{LJZ}. For a quadratic $a$ with $a,a'$ multiplicatively independent the
family is $(a^{c_1}a'^{c_2},(c_1+c_2)r_i)$, with translations $c_2\log(a/a')$ at the
identity and $c_1\log(a/a')$ (up to $2\pi i$) at the conjugate embedding. Complex or
negative conjugates cause no difficulty: the analytic estimates are complex, and the
branches are fixed by the reduction modulo $2\pi i$.
\end{remark}

\begin{remark}[Where the exponent $r+2$ comes from]\label{rem:why}
To have few translations at $\tau$, the family of centers must be long in the direction
$(1,\gamma_\tau)$; as $\tau$ varies these directions span $\Pi\otimes\Q$, of dimension at
most $r+1$, so the family has about $Q^{r+1}$ centers on about $Q^r$ curves at each
embedding. The heights force $Q/w_0\to0$, the collisions $w_0A^mN\to0$, and the size
$v_0\to\infty$; these give $A^{r+2}<\theta^{r+1}$, i.e.\ exponent $r+2$. This is a
heuristic about the present construction, not a theorem. Exponent two would require
removing the height condition.
\end{remark}

\end{document}